\documentclass[reqno]{amsart}
\usepackage{setspace}
\usepackage{graphicx}
\usepackage{amsmath}
\usepackage{amsfonts}
\usepackage{amssymb}
\usepackage[greek,english]{babel}
\begin{document}

\centerline{\Huge \bf Borsuk's Problem in Metric Spaces}

\bigskip
\centerline{\Large Jun Wang, Fei Xue and Chuanming Zong}

\vspace{0.6cm}
\centerline{\begin{minipage}{12.5cm}
{\bf Abstract.} In 1933, K. Borsuk proposed the following problem: Can every bounded set in $\mathbb{E}^n$ be divided into $n+1$ subsets of smaller diameters? In 1965, V. G. Boltyanski and I. T. Gohberg made the following conjecture: Every bounded set in an $n$-dimensional metric space can be divided into $2^n$ subsets of smaller diameters. In this paper, we prove the following result: Every bounded set in an $n$-dimensional metric space can be divided into $2^{n}((n+1)\log (n+1)+(n+1)\log \log (n+1)+5n+5)$ subsets of smaller diameters.
\end{minipage}}



\vspace{0.6cm}
\noindent
{\LARGE\bf 1. Introduction}

\bigskip\noindent
Let $\mathbb{E}^{n}$ be the $n$-dimensional Euclidean space. For any bounded set $S\subseteq \mathbb{E}^{n}$, we call
$$d(S)=\sup\{\|\mathbf{x}, \mathbf{y}\|: \ \mathbf{x}, \mathbf{y}\in S\}$$
the {\it diameter} of $S$, where $\|\cdot\|$ is the {\it Euclidean norm}. Let $b(S)$ be the smallest number such that $S$ can be divided into $b(S)$ subsets with diameters strictly smaller than $d(S)$. In 1933, K. Borsuk \cite{Borsuk1933} proved that an $n$-dimensional ball $B^{n}$ in $\mathbb{E}^{n}$ can not be partitioned into $n$ parts of smaller diameters, which was announced in his ICM-Zurich talk (see \cite{Borsuk1932}). Meanwhile, he proved that every bounded set in $\mathbb{E}^{2}$ can be divided into three subsets of smaller diameters. Based on this fact, he proposed the following problem (Usually, the positive assertion of the problem is known as {\it Borsuk's partition conjecture}):

\bigskip
\noindent
{\bf Borsuk's problem}. {\it Is it true that $$b(S)\leq n+1$$ holds for every bounded set $S$ in $\mathbb{E}^{n}$}?

\medskip
In 1934, T. Bonnesen and W. Fenchel \cite{Bonnesen} proved that every bounded set $S$ of $\mathbb{E}^{2}$ can be divided into three subsets $S_1$, $S_2$ and $S_3$ satisfying
$$d(S_i)\le \frac{\sqrt{3}}{2}\cdot d(S),\quad i=1, 2, 3.$$
In fact, the upper bound $\sqrt{3}/2$ can be attained at circular domains.

In 1947, J. Perkal \cite{Perkal} sketched a proof for a positive answer to the three-dimensional partition problem.  Afterwards, different proofs for this case were discovered by H. G. Eggleston \cite{Eggleston1955}, B. Gr\"unbaum \cite{Grun57}, A. Happes \cite{Heppes} and others (see \cite{grun63,Zong1996}). In particular, B. Gr\"unbaum proved that every bounded set $S$ in $\mathbb{E}^{3}$ can be divided into four parts $S_1$, $S_2$, $S_3$ and $S_4$ satisfying
$$d(S_i)\le 0.9887 \cdot d(S),\quad i=1, 2, 3, 4.$$

In 1993, J. Kahn and G. Kalai \cite{Kahn} surprised the mathematical community by discovering counterexamples to Borsuk's conjecture in high dimensions. They proved that there exist sets $S$ in $\mathbb{E}^{n}$ satisfying
$$b(S)\geq 1.07^{\sqrt{n}}.$$
Therefore, the first counterexample to Borsuk's conjecture occurs in $\mathbb{E}^{21801}.$ Afterwards,  Kahn and Kalai's breakthrough was simplified by N. Alon \cite{alon94} and improved by many authors, in particular by A. Hinrichs and C. Richter \cite{hinr03} to $n\ge 298$ in 2003. In 2014, A. Bondarenko \cite{bond14} presented a $65$-dimensional counterexample to Borsuk's conjecture. Soon after, T. Jenrich and A. E. Brouwer \cite{Jenrich} discovered a $64$-dimensional one. Up to now, Borsuk's problem is still open for $4\leq n \leq 63$. Recently, C. Zong \cite{Zong2021} proposed a computer proof program to this problem.

When Borsuk's conjecture is no longer true in high dimensions, to obtain sharp upper bounds for the partition numbers turns out to be important and interesting. In 1955, H. Lenz \cite{Lenz} proved that
$$b(S)\leq (\sqrt{n}+1)^{n}$$
holds for every bounded set $S$ in $\mathbb{E}^{n}$. This bound was successively improved by L. Danzer \cite{Danzer}, M. Lassak \cite{Lassak} and  O. Schramm \cite{Schramm}. The best-known upper bound is
$$b(S)\leq 5n^{\frac{3}{2}}(4+\log n)\left(\frac{3}{2}\right)^{\frac{n}{2}},$$
which was discovered by O. Schramm in 1988.

Let $(\mathbb{R}^{n}, \|\cdot\|^*)$ be an $n$-dimensional {\it metric space}, i.e., an $n$-dimensional real linear space $\mathbb{R}^{n}$ with norm $\|\cdot\|^*$. It is well-known that $B^*=\{\mathbf{x}\in \mathbb{R}^{n}:\ \|\mathbf{x}\|^*\leq 1\}$ is a centrally symmetric convex body centered at the origin. Usually, $B^*$ is called the unit ball of $(\mathbb{R}^{n}, \|\cdot\|^*)$. On the other hand, if $C$ is a centrally symmetric convex body centered at the origin, then
$$\| {\bf x}, {\bf y}\|^*=\min\{r:\ r\ge 0,\ {\bf x}-{\bf y}\in rC\}$$
defines a metric on $\mathbb{R}^{n}$ and produces a metric space. Therefore, there is a one-to-one correspondence between $n$-dimensional metric spaces and $n$-dimensional centrally symmetric convex bodies centered at the origin. So, for convenience, in this paper we use $M_{C}=\{\mathbb{R}^{n}, \|\cdot\|_{C}\}$ to denote an $n$-dimensional metric space which takes $C$ as the unit ball. For basic concepts and results in metric spaces we refer to \cite{Martini}.

Let $S$ be a bounded set in $M_C$. We define
$$d_C(S)=\sup\{\|\mathbf{x},\mathbf{y}\|_C:\ \mathbf{x},\ \mathbf{y}\in S\}$$
to be the diameter of $S$ and define $b_C(S)$ to be the smallest number such that $S$ can be divided into $b_C(S)$ subsets, all of them having diameters strictly smaller than $d_C(S)$.

\medskip
It is natural to study the analogies of Borsuk's problem in metric spaces. For convenience, let $S$ be a bounded set in a metric space $M_C$ and let $\overline{S}$ denote the completion of the {\it convex hull} of $S$. In 1957, B. Gr\"{u}nbaum \cite{Grunbaum} studied Borsuk's partition problem in metric planes. He showed that, for every bounded set $S$ in a metric plane $M_C$,
$$b_C(S)\le 4 ,$$
where the equality holds if and only if $\overline{S}$ and $C$ are homothetic parallelograms. In 1965, V. G. Boltyanski and I. T. Gohberg \cite[p.~75, 92]{Boltyanski65} made the following conjecture:

\bigskip
\noindent
{\bf Boltyanski-Gohberg conjecture.} {\it For every bounded set $S$ in an $n$-dimensional metric space $M_C$, we have
$$b_C(S)\le 2^n,$$
where the equality holds if and only if $\overline{S}$ and $C$ are homothetic parallelopiped.}

\bigskip
\noindent
{\bf Remark 1.1.} When $C$ is an $n$-dimensional cube, it is well-known and easy to see that
$$b_C(C) = 2^n.$$

\medskip
There are some partial results on the Boltyanski-Gohberg conjecture (see \cite{Boltyanski1977,Hujt14,Lang17,Lian,Wang,Yu,Zong2008}). However, it is still open for all $n\ge 3$. In 1997, C. A. Rogers and C. Zong \cite{Rogers} proved that
$$b_C(S)\leq \binom{2n}{n}(n\log n+n\log \log n+5n)$$
holds for every bounded set $S$ in an $n$-dimensional metric space $M_C$. Up to now, this is the best-known general upper bound for $b_C(S)$. In this paper, we will prove the following theorem:

\bigskip
\noindent
{\bf Theorem 1.1.} {\it For every bounded set $S$ in an $n$-dimensional metric space $M_C$, we have}
$$b_C(S)\leq 2^{n}((n+1)\log (n+1)+(n+1)\log \log (n+1)+5n+5).$$

\vspace{0.6cm}
\noindent
{\LARGE\bf 2. Hadwiger's Covering Conjecture}

\bigskip\noindent
Let $K$ denote a {\it convex body}, a compact and convex set with nonempty interior $\mathrm{int}(K)$, in $\mathbb{E}^n$.
Let $h(K)$ denote the smallest number of translates of $\lambda K$ $(0<\lambda< 1)$ (or $\mathrm{int}(K)$) such that their union contains $K$. In 1957, H. Hadwiger \cite{Hadwiger1957} proposed the following conjecture, which has a close relation with the Boltyanski-Gohberg conjecture.

\bigskip
\noindent
{\bf Hadwiger's covering conjecture.} {\it For every $n$-dimensional convex body $K$, we have
$$h(K)\le 2^n,$$
where the equality holds if and only if $K$ is a parallelopiped.}

\bigskip
This conjecture has been studied by many authors, including K. Bezdek, V. G. Boltyanski, M. Lassak, F. W. Levi, H. Martini, C. A. Rogers, V. P. Soltan, S. Wu and C. Zong. The two-dimensional case was solved by F. W. Levi \cite{Levi} in 1954. However, the conjecture is still open for all $n\ge 3$. Since the target of this paper is the Boltyanski-Gohberg conjecture, we will not go to the details of Hadwiger's conjecture. We refer the interested readers to the references of \cite{Bezdek18,Boltyanski97,Zong2010}. Next, we will introduce two results which will be useful in this paper.

\bigskip\noindent
{\bf Lemma 2.1 (V. G. Boltyanski and I. T. Gohberg \cite{Boltyanski65})} {\it For every bounded set $S$ in a metric space $M_C=\{\mathbb{R}^n,\|\cdot \|_C\}$, we have}
$$b_C(S)\leq h(\overline{S}).$$

\bigskip\noindent
{\bf Lemma 2.2 (C. A. Rogers and C. Zong \cite{Rogers}).} {\it For every $n$-dimensional centrally symmetric convex body $C$, we have}
$$ h(C)\leq 2^{n}(n\log n+n\log \log n+5n).$$

\vspace{0.6cm}
\noindent
{\LARGE \bf 3. Proof of the Theorem}

\bigskip\noindent
First, let us recall a basic concept in convex geometry. For every $n$-dimensional convex body $K$ we define
$$D(K)=K-K=\{{\bf x}-{\bf y}:\ {\bf x},\ {\bf y}\in K\}.$$
Usually, it is known as the {\it difference body} of $K$. Clearly, $D(K)$ is centrally symmetric and convex. In fact the metric $\|\cdot \|_{D(K)}$ defined by $D(K)$ plays the key role of our proof.

Furthermore, for convenience, let $\mathcal{K}^{n}$ denote the set of all $n$-dimensional convex bodies $K$ and let $\mathcal{C}^{n}$ denote the set of all $n$-dimensional centrally symmetric convex bodies $C$. It is easy to see that
$$d_C(S)=d_C(\overline{S})$$
and
$$b_C(S)\le b_C(\overline{S})\eqno(3.0)$$
holds for all bounded sets $S$ in $M_C$. Therefore, to study the Boltyanski-Gohberg conjecture, it is sufficient to deal with the convex bodies in $\mathcal{K}^{n}$.

\bigskip\noindent
{\bf Lemma 3.1.} {\it In every $n$-dimensional metric space $M_C$, we have}
$$\max_{K\in\mathcal{K}^{n}}b_C(K)\leq\max_{K\in\mathcal{K}^{n}}b_{D(K)}(K).$$

\medskip\noindent
{\bf Proof.} Without loss of generality, let $K$ be an $n$-dimensional convex body in $M_C$ with $d_C(K)=1$. Assume that ${\bf x}={\bf x}_1-{\bf x}_2$ and ${\bf y}={\bf y}_1-{\bf y}_2$ are two points in $D(K)$, where all ${\bf x}_1$, ${\bf x}_2$, ${\bf y}_1$ and ${\bf y}_2$ are points in $K$. Then, one can deduce that
$$\|{\bf x},{\bf y}\|_C\le \|{\bf x}_1-{\bf x}_2, {\bf o}\|_C+\|{\bf y}_1-{\bf y}_2, {\bf o}\|_C=\|{\bf x}_1,{\bf x}_2\|_C+\|{\bf y}_1,{\bf y}_2\|_C\le 2$$
and therefore
$$d_C(D(K))\le 2\eqno (3.1)$$
(In fact, one has $d_C(D(K))= 2$). Then, for every point ${\bf x}\in D(K)$, we have
$$\|{\bf x}, {\bf o}\|_C=\mbox{${1\over 2}$}\|{\bf x}, -{\bf x}\|_C\le \mbox{${1\over 2}$}d_C(D(K)) \le 1$$
and consequently
$$D(K)\subseteq C.\eqno (3.2)$$

On the other hand, for any pair of points ${\bf x}$, ${\bf y}\in K$, we have
$$\|{\bf x}, {\bf y}\|_{D(K)}=\mbox{${1\over 2}$}\|{\bf x}-{\bf y}, {\bf y}-{\bf x}\|_{D(K)}$$
and therefore
$$d_{D(K)}(K)=1.\eqno(3.3)$$
If $K$ can be divided into $b_{D(K)}(K)$ subsets $X_1$, $X_2$, $\ldots $, $X_{b_{D(K)}(K)}$ satisfying
$$d_{D(K)}(X_i)<1,\quad i=1, 2, \ldots , b_{D(K)}(K),$$
it follows by (3.2) that
$$d_C(X_i)<1,\quad i=1, 2, \ldots , b_{D(K)}(K)$$
and therefore
$$b_C(K)\leq b_{D(K)}(K).$$
Consequently, we get
$$\max_{K\in\mathcal{K}^{n}}{b_C(K)}\leq \max_{K\in\mathcal{K}^{n}}b_{D(K)}(K).$$
Lemma 3.1 is proved. \hfill{$\Box$}

\medskip
Assume that $K$ is a convex body in $\mathbb{R}^n$. We embed it into $\mathbb{R}^{n+1}=\mathbb{R}^n+\mathbb{R}$ and create a centrally symmetric convex body $K^\bullet$ in $\mathbb{R}^{n+1}$. Setting $K$ in the $n$-dimensional hyperplane
$$H=\{(x_1, x_2, \ldots , x_{n+1}):\ x_{n+1}=0\}$$
of $\mathbb{R}^{n+1}$ and writing ${\bf e}=(0, 0, \ldots , 0, 1)$, then we define
$$K^\bullet=\overline{(K+{\bf e})\cup (-K-{\bf e})}.$$
Clearly, $K^\bullet$ is a centrally symmetric convex body in $\mathbb{R}^{n+1}$.

\bigskip\noindent
{\bf Lemma 3.2.}
$$b_{D(K^\bullet)}(K^\bullet)= b_{D(K^\bullet)}(K+{\bf e})+b_{D(K^\bullet)}(-K-{\bf e})=2b_{D(K^\bullet)}(K+{\bf e}).$$

\medskip\noindent
{\bf Proof.} First of all, we note that
$$d_{D(K^\bullet)}(K^\bullet)=1$$
and
$$\|\mathbf{x},\mathbf{y}\|_{D(K^\bullet)}=1,$$
whenever $\mathbf{x}\in K+{\bf e}$ and $\mathbf{y}\in -K-{\bf e}$. Thus, if $X_1$, $X_2$, $\ldots$, $X_{b_{D(K^\bullet)}(K^\bullet)}$ is a partition of $K^\bullet$ such that
$$d_{D(K^\bullet)}(X_i)<1,\quad i=1,2,\ldots, b_{D(K^\bullet)}(K^\bullet),$$
none of the parts can contain two points $\mathbf{x}\in K+{\bf e}$ and $\mathbf{y}\in -K-{\bf e}$ simultaneously. Therefore, we have $$b_{D(K^\bullet)}(K^\bullet)\geq b_{D(K^\bullet)}(K+{\bf e})+b_{D(K^\bullet)}(-K-{\bf e}).\eqno(3.4)$$

It is well-known in convex geometry that
$$D(K^\bullet)\cap H=D(K).$$
Therefore, we have
$$d_{D(K^\bullet)}(K+{\bf e})=d_{D(K^\bullet)}(-K-{\bf e})=1.\eqno(3.5)$$
Assume that $b_{D(K^\bullet)}(K+{\bf e})=m$ and $K_1$, $K_2$, $\ldots$, $K_m$ is a partition of $K+{\bf e}$ satisfying
$$d_{D(K^\bullet)}(K_i)<1, \quad i=1, 2, \ldots, m.$$
Clearly, $-K_1$, $-K_2$, $\ldots$, $-K_m$ is a partition of $-K-{\bf e}$ satisfying
$$d_{D(K^\bullet)}(-K_i)<1, \quad i=1, 2, \ldots, m$$
and therefore
$$b_{D(K^\bullet)}(-K-{\bf e})=m=b_{D(K^\bullet)}(K+{\bf e}).\eqno(3.6)$$

We define
$$H^+=\{(x_{1},...,x_{n+1}):x_{n+1}\geq 0\},$$
$$H^-=\{(x_{1},...,x_{n+1}):x_{n+1}\leq 0\},$$
$$T_{i}=\overline{(K_{i}\cup (-K-{\bf e}))}\cap H^+$$
for $i=1,$ $2,$ $\ldots,$ $m$, and
$$T_{m+i}=\overline{((-K_{i})\cup (K+{\bf e}))}\cap H^-$$
for $i=1,$ $2,$ $\ldots,$ $m$. It is obvious that
$$K^\bullet=\bigcup_{i=1}^{2m}T_{i}.\eqno(3.7)$$

Next, we proceed to verify that
$$d_{D(K^\bullet)}(T_{i})<1$$
holds for all $i=1,$ $2,$ $\ldots$, $m$. For every pair of points $\mathbf{x}, \mathbf{y}\in T_{i}$, there exist two numbers $\lambda,\ \mu\in [0,\frac{1}{2}]$ and four points $\mathbf{x}_{1}, \mathbf{y}_{1}\in K_{i}$, $\mathbf{x}_{2}, \mathbf{y}_{2}\in -K-{\bf e}$ such that $$\mathbf{x}=(1-\lambda)\mathbf{\mathbf{x}_{1}}+\lambda\mathbf{x}_{2}$$
and
$$\mathbf{y}=(1-\mu)\mathbf{\mathbf{y}_{1}}+\mu\mathbf{y}_{2}.$$
Hence, we have
$$\|\mathbf{x}_{1}-\mathbf{y}_{1}\|_{D(K^\bullet)}<1,$$
$$\|\mathbf{y}_{1}-\mathbf{x}_{2}\|_{D(K^\bullet)}=1$$
and
$$\|\mathbf{x}_{2}-\mathbf{y}_{2}\|_{D(K^\bullet)}\leq 1.$$
Without loss of generality, we assume that $\lambda>\mu$. Then, we get
\begin{align*}
\|\mathbf{x}-\mathbf{y}\|_{D(K^\bullet)}
& =\|(1-\lambda)(\mathbf{x}_{1}-\mathbf{y}_{1})+(\mu-\lambda)(\mathbf{y}_{1}-\mathbf{x}_{2})+\mu(\mathbf{x}_{2}-\mathbf{y}_{2})\|_{D(K^\bullet)}\\
&\leq (1-\lambda)\|\mathbf{x}_{1}-\mathbf{y}_{1}\|_{D(K^\bullet)}+(\lambda-\mu)\|\mathbf{y}_{1}-\mathbf{x}_{2}\|_{D(K^\bullet)}+
\mu\|\mathbf{x}_{2}-\mathbf{y}_{2}\|_{D(K^\bullet)}\\
&<1-\lambda +\lambda -\mu +\mu =1
\end{align*}
and therefore
$$d_{D(K^\bullet)}(T_{i})<1, \quad i=1, 2, \ldots, m.\eqno(3.8)$$
Similarly, one can deduce that
$$d_{D(K^\bullet)}(T_{i})<1, \quad i=m+1, m+2, \ldots, 2m.\eqno(3.9)$$

As a conclusion of (3.6), (3.7), (3.8) and (3.9), we get
$$b_{D(K^\bullet)}(K^\bullet)\leq 2m =b_{D(K^\bullet)}(K+{\bf e})+b_{D(K^\bullet)}(-K-{\bf e}).\eqno(3.10)$$

By (3.4) and (3.10), Lemma 3.2 is proved. \hfill{$\Box$}

\bigskip\noindent
{\bf Lemma 3.3.}
$$b_{D(K^\bullet)}(K+{\bf e})=b_{D(K)}(K).$$

\medskip\noindent
{\bf Proof.} We recall that
$$H=\{(x_{1},...,x_{n+1}):\ x_{n+1}=0\}$$
is an $n$-dimensional hyperplane in $\mathbb{R}^{n+1}$. It is well-known in convex geometry that
$$D(K^\bullet)\cap H=D(K).$$
Thus, when we measure the diameters of subsets of $K+{\bf e}$ by $\| \cdot \|_{D(K^\bullet)}$, the real metric is $\|\cdot \|_{D(K)}$.
Therefore, we get
$$b_{D(K^\bullet)}(K+{\bf e})=b_{D(K)}(K).$$
Lemma 3.3 is proved. \hfill{$\Box$}

\bigskip\noindent
{\bf Proof of Theorem 1.1.} By Lemma 3.2 and Lemma 3.3, we have
$$b_{D(K^\bullet)}(K^\bullet)=2b_{D(K^\bullet)}(K+{\bf e})=2b_{D(K)}(K).\eqno(3.11)$$
Since $K^\bullet\in\mathcal{C}^{n+1}$, by Lemma 2.1 and Lemma 2.2, we get
$$b_{D(K^\bullet)}(K^\bullet)\leq h(K^\bullet)\leq 2^{n+1}((n+1)\log (n+1)+(n+1)\log \log (n+1)+5n+5).$$
Therefore, by (3.11) we have
$$b_{D(K)}(K)=\frac{1}{2}b_{D(K^\bullet)}(K^\bullet)\leq 2^{n}((n+1)\log (n+1)+(n+1)\log \log (n+1)+5n+5)$$
for all $K\in\mathcal{K}^{n}$. Then, Theorem 1.1 follows from (3.0) and Lemma 3.1. \hfill{$\Box$}

\vspace{0.6cm}\noindent
{\bf Acknowledgement.} The work of J. Wang and C. Zong is supported by the National Nature Science Foundation of China (NSFC 11921001) and the National Key Research and Development Program of China (2018YFA0704701). The work of F. Xue is supported by the National Nature Science Foundation of China (NSFC 12201307) and the Natural Science Foundation of Jiangsu Province (BK20210555). J. Wang is the first author. Both F. Xue and C. Zong are corresponding authors.

\vspace{0.6cm}\noindent
Jun Wang, Center for Applied Mathematics, Tianjin University, Tianjin, P. R. China 300072.

\noindent
Email: kingjunjun@tju.edu.cn

\medskip\noindent
Fei Xue, School of Mathematical Sciences, Nanjing Normal University, Nanjing, P. R. China 210046.

\noindent
Email: 05429@njnu.edu.cn

\medskip\noindent
Chuanming Zong, Center for Applied Mathematics, Tianjin University, Tianjin, P. R. China 300072.

\noindent
Email: cmzong@tju.edu.cn

\end{document}